\documentclass[a4paper]{article}
\usepackage{verbatim}
\usepackage{amssymb}
\usepackage{amsmath}
\usepackage{amsthm}
\usepackage{amsfonts}
\usepackage{url}
\usepackage{enumerate}

\usepackage{diagrams} 

\usepackage{mt} 

\usepackage{indent}
\newenvironment{subproof}[1][\proofname]{
  \begin{indentation}{1em}{1em}
    \begin{proof}[#1]
      }{%
        \end{proof}
        \end{indentation}
      }

\theoremstyle{plain}
\newtheorem{thm}{Theorem}[section]

\newtheorem{lemma}[thm]{Lemma}
\newtheorem{prop}[thm]{Proposition}

\newtheorem{fact}[thm]{Fact}
\newtheorem*{fact*}{Fact}
\newtheorem{cor}{Corollary}[thm]

\newtheorem{claim}{Claim}[thm]
\newtheorem*{claim*}{Claim}

\theoremstyle{definition}
\newtheorem{defn}{Definition}[section]
\newtheorem*{defn*}{Definition}

\newtheorem*{definition*}{Definition}

\newtheorem*{defns*}{Definitions}

\newtheorem*{definitions*}{Definitions}
\newtheorem{notation}[defn]{Notation}
\newtheorem*{notation*}{Notation}

\theoremstyle{remark}
\newtheorem{remark}{Remark}[section]
\newtheorem*{remark*}{Remark}

\newtheorem*{example*}{Example}

\newtheorem*{remarks*}{Remarks}

\newtheorem*{examples*}{Examples}

\newtheorem*{note*}{Note}

\author{Martin Bays and Boris Zilber}
\title{Covers of Multiplicative Groups of Algebraically Closed Fields
  of Arbitrary Characteristic}

\date{4 Jan 2011}

\usepackage{covers}

\usepackage[all,cmtip]{xy}

\theoremstyle{plain}
\newtheorem*{prop*}{Proposition}
\newtheorem*{thm*}{Theorem}

\numberwithin{equation}{section}

\begin{document}

\maketitle

\bibliographystyle{alpha}

\begin{abstract}
  We show that algebraic analogues of universal group covers, surjective group
  homomorphisms from a $\mathbb{Q}$-vector space to $F^{\times}$ with
  ``standard kernel'', are determined up to isomorphism of the algebraic
  structure by the characteristic and transcendence degree of $F$ and, in
  positive characteristic, the restriction of the cover to finite fields.
  This extends the main result of ``Covers of the Multiplicative Group of an
  Algebraically Closed Field of Characteristic Zero'' (B. Zilber, JLMS 2007),
  and our proof fills a hole in the proof given there.
\end{abstract}

\section{Introduction}
This paper was conceived as an extension of the main results of
\cite{ZCovers} to fields of positive characteristic. But in the course of
proving the main result a gap in the proof of the main technical theorem,
Theorem 2 (in the case of $n>1$ fields), was detected. So the aim of this
paper has become twofold -- to fix the proof in the characteristic zero case
and to extend it to all characteristics. This goal has now been achieved.

The reader will see that we had to correct the formulation of the theorem of
\cite{ZCovers}. Theorem \ref{thm:determined} below now requires that the
fields $L_1,\dots,L_n$ are from an \defnstyle{independent system}, in the same
sense as in \cite[Section 4]{ZCovers2}, and in accordance with Shelah's theory
of excellence. Indeed, the necessity of this condition has stressed again the
amazingly tight interaction of field-theoretic algebra and very abstract model
theory.

A simple but instructive case of Theorem \ref{thm:determined} is the following
statement:

Let $L_1$ and $L_2$ be linearly disjoint algebraically closed
subfields of a common field of characteristic zero and $L_1L_2$ their
composite. Then the multiplicative group $\multgrpof{(L_1L_2)}$ of the
composite is of the form $A \times (\multgrpof{L_1} \cdot \multgrpof{L_2})$, for
some locally free Abelian group $A$. Surprisingly, even this was apparently
unknown. 

In characteristic $p$ the statement is true with $A$ a locally free
$\Ints[\frac{1}{p}]$-module written multiplicatively.

Here \defnstyle{locally free}, also known as \defnstyle{$\aleph_1$-free},
means that any finite rank subgroup (submodule) is free as an Abelian group
(module). Note that this definition does not agree with the definition of
``locally free'' in general group theory.

Our main technical proposition, Proposition \ref{propn:main}, exhibits a
construction which produces fields $K$ with the multiplicative group of the
form $A\times D$, where $A$ is locally free and $D$ posesses $n$-roots of
elements, for any $n$.  This construction is suggested by Shelah's notion of
independent system and plays a crucial role in proving the uniqueness of
universal covers of the multiplicative group of an algebraically closed field.

\section{Statement of results and outline of proof}
The main theorem of \cite{ZCovers} is:

\begin{thm}\label{thm:categorical0}
  For each cardinal $\kappa$ there is up to isomorphism a unique
  2-sorted structure $\seq{\seq{V;+};\seq{F;+,*};\ex:V\maps F}$ with $V$ a
  divisible torsion-free Abelian group and $F$ an
  algebraically closed field of transcendence degree $\kappa$ such that
  \begin{equation}
    \label{exactSeqChar0}
    \xymatrix{ 0 \ar[r] &\Ints \ar[r] &V \ar[r]^{\ex} &\Fmult \ar[r] &1}
  \end{equation}
  is an exact sequence of groups.
\end{thm}

In positive characteristic the statement must be modified:

\begin{thm}\label{thm:categorical}
  Given a choice of structure $\mathfrak{C_0} :=
  \seq{\seq{\Rats;+};\Fpalg;\ex_0:\Rats\maps\mu}$, where $\mu =
  \multgrpof{(\Fpalg)}$, such that
  \begin{equation}\label{eqn:ex_0}
    \xymatrix{ 0 \ar[r] &\Zpd \ar[r] &\Rats \ar[r]^{\ex_0} &\mu \ar[r] &1} ,
  \end{equation}
  is an exact sequence of groups, for each cardinal $\kappa$ there is
  up to isomorphism a unique 2-sorted structure
  $\mathfrak{C} := \seq{\seq{V;+};\seq{F;+,*};\ex:V\maps F}$ extending
  $\mathfrak{C_0}$ with $V$ a divisible torsion-free Abelian group and $F$ an
  algebraically closed field of characteristic $p$ and transcendence degree
  $\kappa$ such that
  \begin{equation}
    \label{exactSeqCharp}
    \xymatrix{ 0 \ar[r] &\Zpd \ar[r] &V \ar[r]^{\ex} &\Fmult \ar[r] &1}
  \end{equation}
  is an exact sequence of groups.
\end{thm}


Theorems \ref{thm:categorical0} and \ref{thm:categorical} are proven by
showing quasiminimal excellence (\cite{ZQMExc}) of the class of models of an
appropriate $L_{\omega_1,\omega}$-sentence, expressing that we have such a
sequence and, in positive characteristic, that $\ex$ is as specified on
$\Rats\cdot\ker(\ex)$.

For reference, we give a quick outline of the main stages in the proof
now.

$p$ is zero or prime, and $\univacf$ is an arbitrary fixed algebraically
closed field of characteristic $p$.

We use a version of Shelah's notion of an {\em independent system}:

\begin{defn}\label{defn:fromIndieSys}
  We say algebraically closed subfields
  $L_1,\ldots,L_n$ of $\univacf$ are
  \defnstyle{from an independent system} iff there exist an algebraically
  independent set $B\subseteq \univacf$ and subsets $B_i\subseteq B$ such that
  $B=\bigcup_i B_i$ and $L_i = \aclInUniv(B_i)$.
\end{defn}

In the case $n=2$, this condition reduces to saying that $L_1$ is linearly
disjoint from $L_2$ over $L_1\cap L_2$.

\begin{defn}
  If $\ctup\in \univacfmult$ is a $k$-tuple, a \defnstyle{division system
  below $\ctup$} consists of a system of roots
  $(\ctup^{\frac{1}{n}})_{n\in\Nats}$ such that $\ctup^1=\ctup$ and
  $(\ctup^{\frac{1}{nm}})^n = \ctup^{\frac{1}{m}}$. For a rational
  $q=\frac{m}{n}$, we define $c_i^q := (c_i^{\frac{1}{n}})^m$.  For
  an $l\times k$ rational matrix $M = (q_{i,j})_{i,j} \in
  \operatorname{Mat}_{l,k}(\Rats)$, we define $\ctup^M$ to be the $l$-tuple
  $(\Pi_j c_j^{q_{i,j}})_i$, and define $\ctup^\Rats :=
  (\ctup^M)_{M\in\operatorname{Mat}_{1,k}(\Rats)} \leq \univacfmult$.

  If $K\leq \univacf$, we say that \defnstyle{division systems below $\ctup$
  are finitely determined over $K$} iff there exists $m\in\Nats$ such that if
  $(\ctup_1^{\frac{1}{n}})_n$ and $(\ctup_2^{\frac{1}{n}})_n$ are division
  systems below $\ctup$ with $\ctup_1^{\frac{1}{m}}=\ctup_2^{\frac{1}{m}}$,
  then for all $n\in\Nats$ we have that $\ctup_1^{\frac{1}{n}}$ and
  $\ctup_2^{\frac{1}{n}}$ have the same field type over $K$.
\end{defn}

We deduce quasiminimal excellence from the following theorem, the analogue of
Theorem 2 of \cite{ZCovers}.

\newcommand{\Lsmult}{\Pi_i \multgrpof{L_i}}

\begin{thm}\label{thm:determined}
  Let $n\geq 1$ and let $L_1,\ldots,L_n$ be algebraically closed subfields of
  $\univacf$ from an independent system. Let $(\atup,\btup) \in \univacfmult$
  be multiplicatively independent over the product $\Lsmult$. Let
  $(\atup^{\frac{1}{n}})_{n\in\Nats}$ be a division system below $\atup$.

  Then division systems below $\btup$ are finitely determined over
  $L_1 L_2 \ldots L_n (\atup^\Rats)$.
\end{thm}

Theorem \ref{thm:determined} will in turn follow by Kummer theory from
the following proposition describing the structure of the
multiplicative groups of finitely generated perfect extensions of
composites of algebraically closed fields from an independent system.

By $\emp$ is meant $\Zpd$ if $p>0$ and $\Ints$ if $p=0$.

\begin{prop}\label{propn:main}
  Let $\univacf$ be an algebraically closed field, and let $L_1,\ldots,L_n
  \leq \univacf$ be algebraically closed subfields from an independent
  system, $n\geq1$. Let $K$ be the perfect closure of a
  finitely generated extension $L_1\ldots L_n(\betatup) \leq \univacf$ of
  $L_1\ldots L_n$.

  Then $\quot{\multgrpof{K}}{\Lsmult}$ is a locally free $\emp$-module.
\end{prop}

Although Proposition \ref{propn:main} will suffice along with some results from
\cite{ZCovers} to prove Theorem \ref{thm:categorical}, we state here a natural
extension.

\begin{prop}\label{propn:main_ext}
  In each of the following situations,
  $\quot{\multgrpof{(K^{\operatorname{per}})}}{H}$ is a locally free $\emp$-module,
  where $K^{\operatorname{per}}$ is the perfect closure of $K$:
  \begin{itemize}
    \item $K$ is a finitely generated extension of the prime field and $H$ is
      the torsion group of $\multgrpof{K}$
    \item $K$ is a finitely generated extension of the field generated by the
      group $\mu$ of all roots of unity and $H=\mu$
    \item $K$ is a finitely generated extension of the composite $L_1\ldots
      L_n$ of algebraically closed fields from an independent system and $H =
      \Lsmult$.
  \end{itemize}
  In the first two cases, and in the third if $K$ is countable or $n=1$,
  $\quot{\multgrpof{(K^{\operatorname{per}})}}{H}$ is free.
\end{prop}

\begin{remark}
  Theorem 2 of \cite{ZCovers} claims the statement of Theorem
  \ref{thm:determined} for arbitrary finite dimensional algebraically closed
  fields $L_i$, with no independence assumption. The proof given there was
  flawed, but we have no counter-example to this statement; it would be
  interesting to determine whether it is true.
\end{remark}

\section{Torsion-free $\emp$-modules}

\begin{defn}\mbox{}
  \begin{itemize}
  \item For $p$ a positive prime, let $\emp$ be the subring $\Zpd$ of $\Rats$.
  \item For $p=0$, let $\emp$ be the ring $\Ints$.
  \end{itemize}
\end{defn}

To prove Theorem \ref{thm:determined}, we will need to work with the
multiplicative groups perfect (i.e.\ definably closed) subfields of
$\univacf$. These have the natural structure of $\emp$-modules. $\emp$-modules behave, even for $p>0$,
very much like Abelian groups ($\Ints$-modules), and we borrow
definitions and developments from the theory of Abelian groups.


In this section $M$ will be a torsion-free $\emp$-module written
additively.

Here, and throughout the paper, we use tuple notation. A tuple is a
sequence $\atup=(a_i)_{i\in\lambda}$. All tuples will be finite, i.e.
$\lambda\in\omega$, unless otherwise specified. We write (slightly abusively)
$\atup\in A$ to mean that $\atup$ is a finite tuple such that $a_i\in A$ for
all $i$. Unary functions lift to tuples co-ordinatewise - for example,
if $f:A\maps B$ is a function, and $\atup\in A$, then $f(\atup) =
(f(a_1),\ldots,f(a_n)) \in B$.

The ring $\emp$ is a principal ideal domain with fraction field $\Rats$, so we
have the usual definitions:

\begin{defn}
  \mbox{}
  \begin{enumerate}[(i)]
  \item The \defnstyle{span} $\spanof{A} \leq M$ of $A\subseteq M$ is the
    $\emp$-submodule generated by $A$.
  \item $\btup$ is \defnstyle{independent} over $A \leq M$ iff 
    \[ \forall \ntup\in\emp.~(\Sigma_i n_ib_i \in A \implies \ntup = \tuple{0}) .\]
    $B \subseteq M$ is independent over $A$ iff every finite tuple
    $\btup\in B$ is independent over $A$.
  \item The \defnstyle{rank} $r(A)$ of $A \leq M$ is the cardinality
    of any maximal independent $B \subseteq A$. This is well-defined.
  \item $M$ is \defnstyle{free} of rank $\kappa$ iff it is isomorphic
    to the direct sum of $\kappa$ copies of $\emp$, equivalently
    if it is the span of an independent set (called a
    \defnstyle{basis} of $M$) of cardinality $\kappa$.
  \item $M$ is \defnstyle{locally free} iff any finite rank submodule is free.
  \item $M$ embeds in its \defnstyle{divisible hull} $\divhull{M} := M
    \otimes_{\emp} \Rats$, a $\Rats$-vector-space, and $A \leq M$
    embeds in the subspace $\divhull{A} := A \otimes_{\emp} \Rats$ of
    $M \otimes_{\emp} \Rats$, and the embeddings commute.
    $\emp$-independence agrees with $\Rats$-independence in the
    divisible hull, and $r(A)$ is the vector space dimension of
    $\divhull{A}$.
  \end{enumerate}
\end{defn}

Our aim will be to show that certain $\emp$-modules are locally free. To this
end we develop the notions of purity and simplicity:

\begin{defn}
  \mbox{}
  \begin{enumerate}[(i)]
  \item The $\defnstyle{pure hull}$ of a submodule $A \leq M$
    is
    $\purehull{A}{M} := \{x\in M | \exists n\in\emp\setminus\{0\}.~nx\in A\}$.
  \item A submodule $A \leq M$ is \defnstyle{pure} in $M$ iff
    $\purehull{A}{M} = A$.
  \item A tuple $\atup \in M$ is \defnstyle{simple in $M$} iff $\atup$ is
    independent and $\spanof{\atup}$ is pure in $M$. If $A\leq M$ is a pure
    submodule, then $\atup \in M$ is \defnstyle{simple in $M$ mod $A$} iff
    $\quot{\atup}{A}$ is simple in the torsion-free $\emp$-module
    $\quot{M}{A}$.
  \end{enumerate}
\end{defn}

\begin{remark}\label{remark:quotienttf}
  For $A \leq M$, the quotient $\emp$-module $\quot{M}{A}$ is
  torsion-free iff $A$ is pure in $M$.
\end{remark}

\begin{remark}
  In the next section we will be considering quotients of
  multiplicative groups of perfect fields by divisible subgroups
  containing the torsion. It follows from Remark
  \ref{remark:quotienttf} that such quotients are torsion-free
  $\emp$-modules.
\end{remark}

%

\begin{lemma}\label{lemma:freeExtns}
  Suppose $A,B,C$ are $\emp$-modules and $B$ is an extension of $A$ by $C$:
  \begin{equation}
    \xymatrix{ A \ar@{^{(}->}[r] & B \ar@{>>}[r]^{\phi} & C}
  \end{equation}
  Then
  \begin{enumerate}[(i)]
    \item If $A$ and $C$ are free, then $B$ is free
    \item If $A$ and $C$ are locally free, then $B$ is locally free
  \end{enumerate}
\end{lemma}
\begin{proof}\mbox{}
  \begin{enumerate}[(i)]
    \item
      Say $(\phi(b_i))_{i\in I}$ is a basis for $C$. Then $(b_i)_{i\in I}$ are
      independent, and $B = A\oplus\spanof{(b_i)_{i\in I}}$. So $B$ is the direct
      sum of free modules, hence is free.
    \item
      Let $B'$ be a finite rank submodule of $B$. Then we have the exact
      sequence:
      \begin{equation}
	\xymatrix{ A\cap B' \ar@{^{(}->}[r] & B' \ar@{>>}[r]^{\phi} & \phi(B')} .
      \end{equation}
      But $A\cap B'$ and $\phi(B')$ are both finite rank and hence free; so
      $B'$ is free by (i).
  \end{enumerate}
\end{proof}

The following facts are standard results on modules over principal ideal
domains:

\begin{fact}\label{fact:fgfree}
  Any finitely generated torsion-free $\emp$-module is free.
\end{fact}

\begin{fact}\label{fact:subfreefree}
  Any submodule of a free torsion-free $\emp$-module is free.
\end{fact}

\begin{lemma}\label{lemma:freepure}
  A torsion-free $\emp$-module $M$ is locally free iff for every
  finite independent $\atup\in M$, the pure hull of $\spanof{\atup}$
  in $M$ is free.
\end{lemma}
\begin{proof}
  The forward direction is immediate from the definition of local freeness.
  For the converse, suppose $A\leq M$ is finite rank. Let $\atup\in A$ be a
  maximal independent set. Then $A$ is contained in the pure hull of
  $\spanof{\atup}$, which is free by assumption. So $A$ is free by Fact
  \ref{fact:subfreefree}.
\end{proof}

The next two lemmas reduce the condition of purity of a finitely
generated submodule to an easily checked condition on the divisibility
of points.

\begin{lemma}\label{lemma:puresimplepoints}
  A finitely generated submodule $A \leq M$ is pure in $M$
  iff every $a\in A$ which is simple in $A$ is simple in
  $M$.
\end{lemma}
\begin{proof}
  The forward implication is clear. Conversely, suppose $A$ is not pure
  in $M$. Say $\alpha \in M \setminus A$, and $m\alpha = a \in A$ for some
  $m\in\emp\setminus\{0\}$. By Fact \ref{fact:fgfree}, $A$ is free, so the
  pure hull of $a$ in $A$ is free of rank 1, say generated by $a'$. Then $a'$
  is simple in $A$ but not in $M$.
\end{proof}

\begin{lemma}\label{lemma:simplepoint}
  An element $a \in M$ is not simple in $M$ iff $l\alpha = a$ for some
  $\alpha\in M$ and some prime $l \neq p$.
\end{lemma}
\begin{proof}
  Suppose $a$ is not simple in $M$. Then $m\beta = na$ for some $\beta\in
  M\setminus\spanof{a}$ and some $m,n\in\emp$.
  Multiplying up the equation by a power of $p$ we can take
  $m,n\in\Ints$, and by changing $\beta$ we can then take $m\notin
  p\Ints$. We may assume $gcd(m,n)=1$. So there exist $s,t\in\Ints$ such
  that $sm+tn=1$. Then $m(t\beta + sa) = a$. Finish by taking $l$ to be a
  prime divisor of $m$.
\end{proof}

We will have to deal with the delicate question of when a quotient of a
locally free torsion-free $\emp$-module $M$ by a pure submodule $B$ is locally
free, and more generally when for a finite tuple $\ctup\in M$ independent over
$B$ we have that the pure hull of $\quot{\ctup}{B}$ in $\quot{M}{B}$ is free.
Note that if $B$ is finitely generated (equivalently, finite rank) then
$\quot{M}{B}$ is locally free, but that the quotient by an infinite rank
submodule need not be locally free.

The following lemma shows that if, in a certain sense, all the ``extra
divisibility'' of $\ctup$ introduced by quotienting by $B$ is explained by a
finite rank portion of $B$, then the pure hull of $\quot{\ctup}{B}$ is indeed
free.

For $D$ an $\emp$-module and $m\in\emp$, we say that $d \in D$ is
``$m$-divisible in $D$'' iff $\exists d'\in D.~md'=d$.


\begin{lemma}\label{lemma:quotfree}
  Let $M$ be a locally free torsion-free $\emp$-module.
  
  Suppose that $A \leq B \leq M$, that $B$ is pure in $M$, and that $A$ is
  finitely generated.

  Let $\ctup \in M$ be independent over $B$.

  Suppose it holds for all $c \in \spanof{\ctup}$ and all $m\in\emp$
  that if $\quot{c}{B}$ is $m$-divisible in $\quot{M}{B}$, then
  already $\quot{c}{A}$ is $m$-divisible in $\quot{M}{A}$.

  Then the pure hull of $\spanof{\quot{\ctup}{B}}$ is free.
\end{lemma}
\begin{proof}
  \newcommand{\etup}{\tuple{e}}
  Say $A = \spanof{\atup}$. By local freeness, the pure
  hull of $\spanof{\atup\ctup}$ is free, say freely generated by
  $\etup$.

  \begin{claim} $\spanof{\quot{\etup}{B}}$ is the pure hull of
    $\spanof{\quot{\ctup}{B}}$ in $\quot{M}{B}$.
  \end{claim}
  \begin{subproof}\mbox{}
    \begin{itemize}
      \item $\spanof{\quot{\etup}{B}} \leq \purehull{\spanof{\quot{\ctup}{B}}}{\quot{M}{B}}$:
      Indeed let $\quot{e}{B} \in \spanof{\quot{\etup}{B}}$. Without loss of
      generality,
      $e\in \spanof{\etup}$. Then $e$ is in the pure hull of
      $\spanof{\atup\ctup}$, so say $s\cdot e = a+c$, where
      $a\in\spanof{\atup}$ and $c\in\spanof{\ctup}$. But $A\leq B$, so
      $s\cdot\quot{e}{B} = \quot{c}{B} \in \spanof{\quot{\ctup}{B}}$.
    \item $\spanof{\quot{\etup}{B}}$ is pure in $\quot{M}{B}$: Indeed
      suppose $m\cdot\quot{\alpha}{B} = \quot{e}{B}$, where
      $e\in\spanof{\etup}$ and $\alpha \in M$. As above, let
      $c\in\spanof{\ctup}$ and $s\in\Ints$ be such that $s\cdot \quot{e}{B} =
      \quot{c}{B}$. Then $sm\cdot\quot{\alpha}{B} = \quot{c}{B}$, so
      by the assumption, for some $\alpha' \in M$ we have
      $sm\cdot\quot{\alpha'}{A} = \quot{c}{A}$. So $sm\cdot\alpha' =
      c+a$ say, so $\alpha' \in \spanof{\etup}$. But
      $sm\cdot(\alpha-\alpha')\in B$, so by purity of $B$, we have
      $\quot{\alpha}{B} = \quot{\alpha'}{B}$. So $\quot{\alpha}{B} \in
      \spanof{\quot{\etup}{B}}$ as required.
    \end{itemize}
  \end{subproof}
  Now $\spanof{\quot{\etup}{B}}$ is finitely generated, and so by Fact
  \ref{fact:fgfree} is free. So the result follows from the Claim.
\end{proof}

\section{Proof of Proposition \ref{propn:main}}

\begin{notation}
  For subfields $F,F'$ of an algebraically closed field $\univacf$, we write $F\join F'$
  for the perfect closure of the compositum in $\univacf$ of $F$ and $F'$ (so
  in model theoretic terms, $F\join F' = \dcl(F\union F')$), and we write
  $F\join\atup$ for $F\join F'$ where $F'$ is the subfield of $\univacf$
  generated over the prime field by $\{a_1,...,a_n\}$.

  $\mu$ refers to the multiplicative group of all roots of unity.
\end{notation}

We make use of some notions from valuation theory. We consider a
\defnstyle{place} of a field $\pi : K \maps k$ to be a partially defined ring
homomorphism such that the domain of definition $\valring_\pi := \dom(\pi)$ is
a valuation ring. If $k\leq K$, we write $\pi : K \maps_k k$ to indicate that
$\pi$ is the identity on $k$ - in other words, that the field embedding of $k$
in $K$ is a section of $\pi$. Such a $\pi$ is sometimes called a
\defnstyle{specialisation} of $K$ to $k$.

We make use of the Newton-Puiseux theorem, or rather the following
generalisation to arbitrary characteristic:

\newcommand{\genPuiseux}[2]{#1\{\{#2\}\}}
\newcommand{\genPowerSeriesQ}[2]{#1((#2^{\Rats}))}
\begin{fact}[Rayner \cite{RaynerNewtonPosChar}, cited in
  \cite{KedlayaNewtonPosChar}]\label{fact:Puiseux}
  Let $L$ be an algebraically closed field of characteristic $p\geq 0$. Let
  $\genPowerSeriesQ{L}{t}$ be the field of generalised formal power series in $t$
  with coefficients in $L$ and rational exponents, and let $\genPuiseux{L}{t}
  \leq \genPowerSeriesQ{L}{t}$ be the subfield consisting of those power
  series with support $S\subseteq \Rats$ satisfying:
  \begin{itemize}
    \item there exists $m\in\Ints\setminus\{0\}$ such that $mS \subseteq \emp$.
  \end{itemize}

  Then $\genPuiseux{L}{t}$ is an algebraically closed field.
\end{fact}

\begin{lemma}\label{lemma:placesAndPreservation}
  Let $L$ be an algebraically closed subfield of an algebraically closed field
  $\univacf$; suppose $L$ contains algebraically closed subfields
  $k_i$ for $i\in\{1,\ldots n\}$; let $\lambda\in \univacf$ be transcendental over $L$; let
  $K := \aclInUniv(L(\lambda))\geq L$, and let $k_i' := \aclInUniv(k_i(\lambda))$.
  Further, let $k_0 \leq L$ be a perfect subfield, and let $k_0' := k_0$.

  Then for any place $\pi : K \maps_L L$ such that $\pi(\lambda) \subseteq
  \bigcap_{i>0} k_i$,
  \[ \pi(\bigjoin_{i\geq 0} k_i') = \bigjoin_{i\geq 0} k_i .\]
\end{lemma}
\begin{proof}
  Since replacing $\lambda$ with $\lambda - \pi(\lambda)$ does not alter
  $K$ or $k_i'$, and $\lambda - \pi(\lambda)$ is also transcendental over $L$,
  we may assume that $\pi(\lambda) = 0$.

  Let $\genPuiseux{L}{\lambda}$ be the field of generalised Puiseux
  series, as defined in Fact $\ref{fact:Puiseux}$.
  Let $\pi' : \genPuiseux{L}{\lambda} \maps L$ be the standard power series
  residue map.

  $\pi'$ agrees with $\pi$ on $L(\lambda)$, so by the Conjugation Theorem
  \cite[3.2.15]{PrestelEnglerValuedFields} we may embed $K$ into
  $\genPuiseux{L}{\lambda}$ over $L(\lambda)$ in such a way that $\pi$ agrees with $\pi'$.

  Now for $i>0$, the
  subfield $\genPuiseux{k_i}{\lambda} \leq \genPuiseux{L}{\lambda}$ of power series with coefficients from $k_i$ is algebraically
  closed and contains $k_i(\lambda)$, so contains $k_i'$. Similarly, $k_0' =
  k_0 \leq \genPuiseux{k_0}{\lambda}$.

  Now 
  \begin{align*}
    \pi(\bigjoin_{i\geq 0} k_i')
    &\leq \pi'( \bigjoin_{i\geq 0} ( \genPuiseux{k_i}{\lambda} ) ) \\
    &\leq \pi'( \genPuiseux{( \bigjoin_{i\geq 0} k_i )}{\lambda} ) \\
    &= \bigjoin_{i\geq 0} k_i
  \end{align*}
\end{proof}

\begin{lemma}\label{lemma:placeExists}
  Suppose $L_1,\ldots,L_n\leq\univacf$ are algebraically closed subfields from
  an independent system, witnessed by an independent set
  $B=B_1\cup\ldots\cup B_n$ as in Definition \ref{defn:fromIndieSys}.
  Let $B^0\subseteq B$ and define $B_i^0 := B_i \cap B^0$ and $L_i^0 :=
  \aclInUniv(B_i^0)$. Let $C\subseteq \aclInUniv(B^0)$.

  Then there exists a place $\pi : \aclInUniv(B) \maps_{\aclInUniv(B^0)}
  \aclInUniv(B^0)$ such that $\pi(L_i) = L_i^0$ and $\pi(\bigjoin L_i \join C)
  = \bigjoin L_i^0\join C$ .

  Furthermore, for any finite tuple $\ctup\in \multgrpof{\aclInUniv(B)}$,
  $\pi$ can be chosen such that $\pi(\ctup)\in \multgrpof{\aclInUniv(B^0)}$.
\end{lemma}
\begin{proof}
  Let the possibly infinite tuple $\btup=(b_\alpha)_{\alpha<\lambda}$
  enumerate $B\setminus B^0$.

  For $\beta\leq\lambda$, define $B^\beta :=
  B^0\union\{b_\alpha|\alpha<\beta\}$; $L^\beta := \aclInUniv(B^\beta)$;
  $B_i^\beta := B_i \cap B^\beta$; $L_i^\beta := \aclInUniv(B_i^\beta)$, and
  $K^\beta := \bigjoin L_i^\beta \join C$.

  Let $f_{i,j}(\btup) \in L^0[\btup]$ be the non-zero coefficients
  of a minimal polynomial in $L^0[\btup][X]$ for $c_i$ over $L^0(\btup)$. Let
  $\atup=(a_\alpha)_{\alpha<\lambda}\in \aclInUniv(\emptyset)$ be such that
  $f_{i,j}(\atup) \neq 0$ for all $i,j$.

  We define, by transfinite recursion on $\beta\leq\lambda$,
  places $\pi^\beta : L^\beta \maps_{L^0}
  L^0$ such that
  $\pi^\beta(b_\alpha)=a_\alpha$ for $\alpha<\beta$, and
  $\pi^\beta(L_i^\beta) = L_i^0$ and
  $\pi^\beta(K^\beta)=K^0$, and
  $\pi^\beta\restriction L^\gamma = \pi^\gamma$ for $\gamma\leq \beta$.

  Define $\pi^0 := \id_{L^0}$, and take unions at limit ordinals. If
  $\beta=\gamma+1$ is a successor ordinal, by Lemma
  \ref{lemma:placesAndPreservation} 
  if
  $\pi^{\gamma+1}_\gamma : L^{\gamma+1} \maps_{L^\gamma} L^\gamma$
  is a place such that
  $\pi^{\gamma+1}_\gamma(b_\gamma) = a_\gamma$, then
  $\pi^{\gamma+1}_\gamma(K^{\gamma+1}) = K^\gamma$; clearly we also have
  $\pi^{\gamma+1}_\gamma(L_i^{\gamma+1}) = L_i^\gamma$.
  So $\pi^{\gamma+1} := \pi^{\gamma}\composition\pi^{\gamma+1}_\gamma$ is as
  required.

  Now let $\pi := \pi^\lambda$. By the condition on $\atup$, we have
  $\pi(c_i)\in\multgrpof{L^0}$.
\end{proof}

\begin{lemma}\label{lemma:specFinGen}
  Let $K\geq L$ be algebraically closed fields, and let $\pi: K \maps_L L$ be a
  place. Let $k_0 \leq K$ be a perfect subfield such that $\pi k_0\leq k_0$.
  Let $k_1\geq k_0$ be a finite extension.
  
  Then there exists a finite extension $k'\geq k_1$ such that $\pi k'\leq
  k'$.
\end{lemma}
\begin{proof}
  We may assume that $k_1/k_0$ is Galois.

  For $i\geq 1$, define $k_{i+1} := k_i ( \pi k_i )$.

  A finite extension of a perfect field is perfect, so each $k_i$, and hence
  each $\pi k_i$, is perfect.

  Normality of a finite field extension implies
  \cite[3.2.16(2)]{PrestelEnglerValuedFields} normality of the corresponding
  extension of residue fields; it follows inductively that for all $i\geq 0$,
  the extensions $k_{i+1}/k_i$ and $\pi k_{i+1} /\pi k_i $ are Galois.

  Now $k_{i+2}$ is generated over $k_{i+1}$ by $\pi k_{i+1}$, and $\pi k_i
  \leq k_{i+1}$, so $[ k_{i+2} : k_{i+1} ] \leq [ \pi k_{i+1} : \pi k_i ]$.
  Also, $[ \pi k_{i+1} : \pi k_i ] \leq [ k_{i+1} : k_i ]$.
  So after some $n$, the degrees reach their minimum level, say
  \[ d = [ \pi k_{n+2} : \pi k_{n+1} ] = [ k_{n+2} : k_{n+1} ]
  = [ \pi k_{n+1} : \pi k_n ] = [ k_{n+1} : k_{n} ] .\]
  By the fundamental inequality of valuation theory
  \cite[3.3.4]{PrestelEnglerValuedFields}, 
  \begin{enumerate}[(I)]
    \item any $\sigma\in\Gal( k_{n+1} / k_n)$ preserves $\valring_\pi \cap
      k_{n+1}$;
    \item any $\sigma\in\Gal( k_{n+2} / k_{n+1})$ preserves $\valring_\pi \cap k_{n+2}$.
  \end{enumerate}

  Now $\pi k_{n+1} = (\pi k_n)(\pi \beta)$ say, some $\beta\in k_{n+1}$.
  Let $\beta=\beta_1,\beta_2,\ldots,\beta_s$ be the $k_n$-conjugates of $\beta$.
  By (I), $\beta_i \in \valring_\pi$ for all $i$. Applying $\pi$ to the
  minimum polynomial $\Pi_i (x-\beta_i)$, we see that $s=d$ and the $(\pi
  k_n)$-conjugates of $\pi \beta$ are precisely $(\pi \beta_i)_i$.

  Now suppose for a contradiction that $\sigma \in \Gal( k_{n+2} / k_{n+1} )
  \setminus \{\id\}$. We have $k_{n+2} = k_{n+1}(\pi \beta)$, so $\sigma(\pi \beta)
  = \pi \beta_i$ some $i>1$.

  Now $\beta - \pi\beta \in \maxideal_\pi \cap k_{n+1}$, but $\sigma(\beta -
  \pi \beta) = \beta - \sigma\pi\beta = \beta - \pi\beta_i \notin \maxideal_\pi
  \cap k_{n+1}$. This contradicts (II).

  So $d=1$, and so $\pi k_n \leq k_n$.
\end{proof}

\begin{fact}\cite[Proposition 1]{MayMultGrps}\label{fact:fgextnacffree}
  Let $E\geq F$ be a finitely generated regular extension. Then
  $\quot{\multgrpof{E}}{\multgrpof{F}}$ is free as an Abelian group.
\end{fact}

This fact slightly extends the second statement of \cite[Lemma 2.1]{ZCovers}.
The proof involves considering the Weil divisors of a normal projective
variety over $F$ with function field $E$.

We translate this result to our context of perfect fields and $\emp$-modules:

\begin{cor}\label{cor:fgextnacffree}
  Let $\percl{E}$ be the perfect closure of a finitely generated regular extension
  $E$ of a perfect field $F$. Then
  $\quot{\multgrpof{\percl{E}}}{\multgrpof{F}}$ is free as an $\emp$-module.
\end{cor}
\begin{proof}
  This is immediate from Fact \ref{fact:fgextnacffree}, on noting that if
  $(\quot{e_i}{\multgrpof{F}})_{i<\kappa}$ is a basis for
  $\quot{\multgrpof{E}}{\multgrpof{F}}$ as an Abelian group, then
  $(\quot{e_i}{\multgrpof{F}})_{i<\kappa}$ is a basis for
  $\quot{\multgrpof{\percl{E}}}{\multgrpof{F}}$ as an $\emp$-module.
\end{proof}

\begin{prop*}[\ref{propn:main}]
  Let $\univacf$ be an algebraically closed field, and let $L_1,\ldots,L_n
  \leq \univacf$ be algebraically closed subfields from an independent
  system, $n\geq1$. Let $\betatup\in \univacf$ be an arbitrary finite tuple,
  and let $K := L_1 \join \ldots \join L_n \join \betatup \leq \univacf$.

  Then $\quot{\multgrpof{K}}{\Lsmult}$ is a locally free $\emp$-module.
\end{prop*}
\begin{proof}
  The $n=1$ case of the proposition follows from Corollary
  \ref{cor:fgextnacffree}; we proceed to prove the proposition by induction on
  $n$.

  Let $B, B_i$ be as in Definition \ref{defn:fromIndieSys}.

  Let $L := L_1$, let $P := \bigjoin_{i>1} L_i$, and let $H := \Pi_{i>1}
  \multgrpof{L_i} \leq \multgrpof{P}$.

  \newcommand{\HL}{H\multgrpof{L}}
  \newcommand{\PLb}{P\join L\join\betatup}
  \newcommand{\Pab}{P\join\atup\join\betatup}
  \newcommand{\kL}{k\join L}
  \newcommand{\PLbmult}{\multgrpof{(\PLb)}}
  \newcommand{\kLmult}{\multgrpof{(\kL)}}
  \newcommand{\kmult}{\multgrpof{k}}

  We first show that we may reduce to the case that $\betatup$ is
  algebraic over $P\join L = \bigjoin_i L_i$.
  Indeed, the relative algebraic closure of $P\join L$ in $P\join L\join
  \betatup$, is an algebraic subextension of the finitely generated extension
  $(P\join L)(\betatup)$ of
  $P\join L$, and so is a finite extension $P\join L \join \betatup'$ say, where
  $\betatup' \in \aclInUniv(P\join L)$.

  By Corollary
  \ref{cor:fgextnacffree}, $\quot{\multgrpof{(P\join L \join
  \betatup)}}{\multgrpof{(P\join L \join \betatup')}}$ is free. So
  by Lemma \ref{lemma:freeExtns}, we need only show that
  $\quot{\multgrpof{(P\join L \join \betatup')}}{\HL}$ is locally free.

  So we suppose that $\betatup \in \aclInUniv(P\join L)$.

  We claim further that we may assume $B$ to be finite.  Indeed, suppose
  $B^0\subseteq_{\operatorname{fin}} B$ is such that
  $\betatup\in\aclInUniv(B^0)$. Let $B_i^0 := B_i \cap B^0$, and define $L_i^0
  := \aclInUniv(B_i^0)$ and $K^0 := \bigjoin L_i^0\join\betatup \leq K$.
  
  Note that $\Lsmult\cap\multgrpof{K^0} = \Pi_i \multgrpof{L_i^0}$. Indeed, if
  $x = \Pi_i a_i \in \Lsmult\cap\multgrpof{K^0}$, then by Lemma
  \ref{lemma:placeExists} there exists a place $\pi_0 : K
  \maps_{K^0} K^0$ such that $\pi_0(a_i)\in \multgrpof{L_i^0}$, so
  $x=\pi_0(x)=\Pi_i \pi_0(a_i) \in \Pi_i\multgrpof{L_i^0}$.

  So the $\emp$-module $M(B^0) := \quot{\multgrpof{K^0}}{\Lsmult}$ is
  isomorphic to $\quot{\multgrpof{K^0}}{\Pi_i\multgrpof{L_i^0}}$.  By the
  existence of $\pi_0$, we have that $M(B^0)$ is pure in
  $M:=\quot{\multgrpof{K}}{\Lsmult}$. So assuming the current lemma for finite
  $B$, we have that $M$ is the union of the locally free pure submodules $M(B^0)$ as $B^0$
  ranges through the finite subsets of $B$ for which
  $\betatup\in\aclInUniv(B^0)$, and so $M$ is locally free as required.

  So we assume that $B$ is finite.

  We aim to apply Lemma
  \ref{lemma:freepure}. So let $\btup \in \PLb$ be multiplicatively
  independent over $\HL$; we want to show that the pure hull of
  $\spanof{\quot{\btup}{\HL}}$ in $\quot{\PLbmult}{\HL}$ is free.

  Let $(c_i)_i$ enumerate $\betatup\btup$.
  \begin{claim}\label{claim:mainindn:placeExists}
    There exist a finitely generated extension $k$ of $P$
    and a place $\pi: \aclInUniv(LP) \maps_L L$
    such that
    \begin{enumerate}[(i)]
      \item $\kL \geq \PLb$;
      \item $\forall i\qsep c_i\in k$;
      \item $L=\acl^L(k\cap L)$;
      \item $\pi(k) = k\cap L$;
      \item $\pi(c_i)\in\multgrpof{L}$.
    \end{enumerate}
  \end{claim}
  \begin{proof}
    By Lemma \ref{lemma:placeExists} with $B^0 := B_1$ and $C:=B_1$, there
    exists a place $\pi : \aclInUniv(LP)\maps_L L$ such that $\pi(P \join B_1)
    = L$ and $\pi(c_i)\in\multgrpof{L}$.

    By Lemma \ref{lemma:specFinGen}, there exists a finite extension $k$ of
    $P\join B_1 \join \ctup$ such that $\pi(k)\leq k$.
    
    Then $k$ and $\pi$ are as required.
  \end{proof}

  \begin{claim}\label{claim:mainindn:cohom}
    If $b\in \kmult$ is simple in $\kmult$ mod
    $(\kmult \cap H\multgrpof{L})$, then $b$ is simple in
    $\kLmult$ mod $H\multgrpof{L}$.

    Furthermore, identifying $\quot{\kmult}{\kmult\cap\HL}$
    with the submodule\\ $\quot{\kmult}{\HL}$ of $\quot{\kLmult}{\HL}$, we
    have that for any $\ctup \in \kmult$ if $\spanof{\quot{\ctup}{\HL}}$
    is pure in $\quot{\kmult}{\HL}$ then it is is pure in
    $\quot{\kLmult}{\HL}$.
  \end{claim}

  \begin{subproof}[Proof of Claim \ref{claim:mainindn:cohom}]
    Suppose $b$ is not simple
    in $\kLmult$ mod $H\multgrpof{L}$. By Lemma
    \ref{lemma:simplepoint} and the fact that $H\multgrpof{L}$ is
    divisible in $\kLmult$, we have $\alpha^q = b$ for
    some $\alpha\in\kLmult \setminus \kmult$ and some
    prime $q\neq p$.

    Now $k(\alpha)$ is a degree $q$ cyclic extension of $k$, so this is a
    Galois extension, $\Galofover{k(\alpha)}{k} \isom \Intsmod{q}$, and
    $k(\alpha)$ is perfect.

    Let $F_0 := k\cap L$ and $F_1 := k(\alpha)\cap L$. Let $F_2 \leq L$ be a
    finite extension of $F_1$ such that $\alpha\in k\join F_2$ and $F_2$ is
    Galois over $F_0$. Note that $F_2\cap k = F_0$ and $F_2\cap k(\alpha) =
    F_1$.

    By \cite[VI Thm 1.12]{LangAlg}, $k\join F_2$ is Galois over $k$ and
    restriction to $F_2$ gives an isomorphism of finite groups

    \[ \restriction_{F_2} : \Galofover{k\join F_2}{k} \maps
    \Galofover{F_2}{F_0} ,\]

    and $\Galofover{F_2}{F_1}$ is the image under $\restriction_{F_2}$ of the
    normal subgroup $\Galofover{k\join F_2}{k(\alpha)}$ of
    $\Galofover{k\join F_2}{k}$.

    So $F_1$ is Galois over $F_0$ and
    \begin{align*}
      \Galofover{F_1}{F_0}
      &\isom \quot{\Galofover{F_2}{F_0}}{\Galofover{F_2}{F_1}}\\
      &\isom \quot{\Galofover{k\join F_2}{k}}{\Galofover{k\join
      F_2}{k(\alpha)}}\\
      &\isom \Galofover{k(\alpha)}{k}\\
      &\isom \Intsmod{q}.
    \end{align*}

    By \cite[VI Thm 1.12]{LangAlg} again, $\Galofover{kF_1}{k} \isom
    \Galofover{F_1}{F_0} \isom \Galofover{k(\alpha)}{k}$. So $k\join F_1 =
    kF_1 = k(\alpha)$, and we have the following lattice diamond:

    \begin{diagram}[small]
      &                       & k(\alpha) &                       &     \\
      & \ruLine^{\Intsmod{q}} &            & \rdLine               &     \\
      k &                       &            &                       & F_1 \\
      & \rdLine               &            & \ruLine_{\Intsmod{q}} &     \\
      &                       & k \cap L  &                       &     \\
    \end{diagram}

    Since the torsion group $\mu$ is contained in
    $\multgrpof{(k\cap L)}$, by \cite[VI 6.2]{LangAlg} $F_1 = (k\cap
    L)(\gamma)$ for some $\gamma$ such that $\gamma^q\in k\cap L$.

    Now $k(\alpha) = k\join F_1 = k(\gamma)$, so say $\gamma =
    \Sigma_{i<q}c_i\alpha^i$, with $c_i\in k$. Let
    $\sigma\in\Galofover{k(\alpha)}{k}$ restrict non-trivially to
    $F_1$. Say $\sigma(\alpha)=\zeta\alpha$ and
    $\sigma(\gamma)=\zeta^l\gamma$, where $(l,q)=1$, and $\zeta$ is a primitive
    $q$th root of unity. So

    \[ \Sigma_{i<q}c_i\zeta^l\alpha^i = \zeta^l\gamma = \sigma(\gamma)
    = \Sigma_{i<q}c_i\zeta^i\alpha^i .\]

    Since $(\alpha^i)_i$ is a basis for the $k$-vector-space
    $k(\alpha)$, we have $\gamma = c_l\alpha^l$.

    Now say $sl + tq = 1$. Then $\gamma^s = c_l^s\alpha b^{-t}$. So
    letting $d := c_l^{-s}b^t \in k$, we have
    \[ d^q\gamma^{sq} = \alpha^q = b .\]

    But $\gamma^{sq} = (\gamma^q)^s \in \multgrpof{(k\cap L)}$, so by Lemma
    \ref{lemma:simplepoint} $b$ is
    not simple in $\multgrpof{k}$ mod $(\multgrpof{k} \cap H\multgrpof{L})$.

    This completes the proof of the first statement. The
    ``Furthermore'' part follows by Lemma
    \ref{lemma:puresimplepoints}.
  \end{subproof}

  We aim to apply Lemma \ref{lemma:quotfree}. Let $N :=
  \quot{\multgrpof{k}}{H}$, which by the inductive hypothesis is a torsion-free locally free
  $\emp$-module; let $D := \quot{(\multgrpof{k}\cap\HL)}{H}$, which is a pure
  submodule of $N$; and let $A := \spanof{\quot{\pi(\btup)}{H}} \leq D$.

  \begin{claim}\label{claim:mainindn:rootsReduce}
    Let $\quot{b}{H} \in \spanof{\quot{\btup}{H}}$, let $m\in\emp$.

    If $\quot{b}{H}$ has an $m$th root modulo $D$ in $N$, then
    $\quot{b}{H}$ has an $m$th root modulo $A$ in $N$.
  \end{claim}
  \begin{subproof}
    Say $\quot{\lambda}{H}(\quot{\alpha}{H})^m = \quot{b}{H}$, where $\alpha \in
    \multgrpof{k}$, and $\lambda\in \HL$. Since $H$ is divisible, we may suppose that $\lambda \in
    \multgrpof{L}$, where $b\in\spanof{\btup}$, and $\lambda\alpha^m = b$.

    Applying $\pi$, we obtain (recalling that $\pi(b_i)\in\multgrpof{L}$
    and that $\pi$ fixes $L\ni\lambda$)
    \[ \lambda = \pi(\lambda) = \quot{\pi(b)}{\pi(\alpha)^m} .\]
    So
    \[ \pi(b)\left(\quot{\alpha}{\pi(\alpha)}\right)^m = b .\]
    But $\pi(b) \in \spanof{\pi(\btup)}$ and
    $\pi(\alpha)\in\pi(k)\subseteq k$, so this shows that
    $\quot{b}{H}$ has an $m$th root modulo $A$ in $N$.
  \end{subproof}
  It follows from Lemma \ref{lemma:quotfree} and Claim
  \ref{claim:mainindn:rootsReduce} that the pure hull of
  $\spanof{\quot{\btup}{\HL}}$ in $\quot{\kmult}{\HL}$ is free;
  by Claim \ref{claim:mainindn:cohom}, the pure hull in $\quot{\PLbmult}{\HL}
  = \quot{\kLmult}{\HL}$ is also free.

  Applying Lemma \ref{lemma:freepure}, this completes the proof of Proposition
  \ref{propn:main}.
\end{proof}

\begin{prop}
  In each of the following situations,
  $\quot{\multgrpof{\percl{K}}}{H}$ is a locally free $\emp$-module:
  \begin{itemize}
    \item $K$ is a finitely generated extension of the prime field and $H$ is
      the torsion group of $\multgrpof{K}$
    \item $K$ is a finitely generated extension of the field generated by the
      group $\mu$ of all roots of unity and $H=\mu$
    \item $K$ is a finitely generated extension of the composite $L_1\ldots
      L_n$ of algebraically closed fields from an independent system and $H =
      \Lsmult$.
  \end{itemize}
  In the first two cases, and in the third if $K$ is countable or $n=1$,
  $\quot{\multgrpof{(K^{\operatorname{per}})}}{H}$ is free.
\end{prop}
\begin{proof}
  In characteristic 0, the first case is the first part of the statement of
  \cite[Lemma 2.1]{ZCovers}, and the second case is \cite[Lemma 2.14(ii)]{ZCovers}.

  In characteristic $p>0$, both the first and second case follow from
  Corollary \ref{cor:fgextnacffree} with $F$ being $K\cap \Fpalg$.

  In all characteristics the third case is precisely Proposition
  \ref{propn:main}. Freeness in the countable case follows from Pontyragin's
  theorem (\cite[19.1]{FuchsInfAb}), and in the $n=1$ case from 
  Fact \ref{fact:fgextnacffree}.
\end{proof}

\section{Proof of Theorem \ref{thm:determined}}
Theorem \ref{thm:determined} will follow from Proposition
\ref{propn:main} by Kummer theory, our use of which is packaged in the
following lemma:

\begin{lemma}\label{lemma:kummer}
  Let $K$ be a perfect field containing the roots of unity $\mu$, and let $F \geq K$
  algebraically closed. Let $\atup \in \multgrpof{K}$ such that
  $\quot{\atup}{\mu}$ is simple in $\quot{\multgrpof{K}}{\mu}$.  Let
  $n\in\Nats$. Then all choices of $\alphatup \in \multgrpof{F}$ such that
  $\alphatup^n = \atup$ have the same field type over $K$.
\end{lemma}
\begin{proof}
  Let $\alphatup$ be such. Say $n = p^tm$ where $(m,p)=1$. Since the field type of $\alphatup$
  is determined by that of $\alphatup^{p^t}$, it suffices to consider the case
  that $t=0$. By Kummer theory (\cite[VI\S 8]{LangAlg}),
  \[ \Galofover{K(\alphatup)}{K} \isom
  \operatorname{Hom}\left(
  \quot{\spanofover{\atup}{\Ints}}{\spanofover{\atup}{\Ints} \cap
  (\multgrpof{K})^n}, \quot{\Ints}{n\Ints}\right)
  \isom
  \quot{\spanofover{\atup}{\Ints}}{\spanofover{\atup}{\Ints} \cap
    (\multgrpof{K})^n} ,\]
  where $(\multgrpof{K})^n$ is the $n$-powers subgroup of $\multgrpof{K}$.

  By simplicity, $\spanofover{\atup}{\Ints} \cap (\multgrpof{K})^n =
  \spanofover{\atup^n}{\Ints}$. So $\Galofover{K(\alphatup)}{K} \isom
  \left(\quot{\Ints}{n\Ints}\right)^{\sizeof{\atup}}$.
\end{proof}
\begin{thm*}[\ref{thm:determined}]
  Let $n\geq 1$ and let $L_1,\ldots,L_n$ be algebraically closed subfields of
  $\univacf$ from an independent system. Let $(\atup,\btup) \in \univacfmult$
  be multiplicatively independent over the product $\Lsmult$. Let
  $(\atup^{\frac{1}{n}})_{n\in\Nats}$ be a division system below $\atup$.

  Then division systems below $\btup$ are finitely determined over
  $L_1 L_2 \ldots L_n (\atup^\Rats)$.
\end{thm*}
\begin{proof}
  Let $\ctup := \atup\btup$.

  Let $K := \bigvee_i L_i \join \ctup$.
  Let $\Gamma_1$ be the pure hull of $\quot{\ctup}{\Lsmult}$ in
  $\quot{\multgrpof{K}}{\Lsmult}$, and let $\Gamma$ be the pure hull of
  $\quot{\ctup}{\mu}$ in $\quot{\multgrpof{K}}{\mu}$.

  Since $\ctup$ is multiplicatively independent over $\Lsmult$ and
  $\mu\leq\Lsmult\leq\multgrpof{K}$, the $\emp$-modules $\Gamma$ and
  $\Gamma_1$ are isomorphic; by Proposition \ref{propn:main} and Lemma
  \ref{lemma:freepure}, $\Gamma_1$, and hence $\Gamma$, is free.

  Now let $m$ be such that $\Gamma^m \leq \spanof{\quot{\ctup}{\mu}}$. Suppose
  $(\btup_1^{\frac{1}{n}})_n$ and $(\btup_2^{\frac{1}{n}})_n$ are division
  systems below $\btup$ such that
  $\btup_1^{\frac{1}{m}}=\btup_2^{\frac{1}{m}}$, and let $k\in\Nats$; we claim
  that $\btup_i^{\frac{1}{k}}$ have the same field type over $\bigvee_i L_i
  \join \atup^\Rats$.

  Define division systems above $\ctup$ by $\ctup_i^{\frac{1}{n}} :=
  \atup^{\frac{1}{n}}\btup_i^{\frac{1}{n}}$. By the choice of $m$, there
  exists $M\in\GL_{\sizeof{\ctup}}(\Rats)$ such that
  $\ctup':=\ctup_1^M=\ctup_2^M$ and $\quot{\ctup_i^M}{\mu}$ is a $\emp$-basis
  for $\Gamma$.

  By Lemma \ref{lemma:kummer}, for all $n\in\Nats$ all choices of
  $\ctup{'}^{\frac{1}{n}}$ have the same field type over $K$.  Hence for any
  $l$, we have that $(\atup^{\frac{1}{l}}\btup_1^{\frac{1}{k}})$ and
  $(\atup^{\frac{1}{l}}\btup_2^{\frac{1}{k}})$ have the same field type over
  $K$, and so $\btup_i^{\frac{1}{k}}$ have the same field type over $K\join
  \atup^{\frac{1}{l}}$. Since this holds for all $l$, we have that
  $\btup_i^{\frac{1}{k}}$ have the same field type over $\bigvee_i L_i \join
  \atup^\Rats$, as required.
\end{proof}

\section{Proof of Theorems \ref{thm:categorical0} and \ref{thm:categorical}}

\newcommand{\qtup}{\tuple{q}}
\newcommand{\CClass}{\mathcal{C}}
\newcommand{\cl}{\operatorname{cl}}
\providecommand{\image}{\operatorname{im}}
\providecommand{\spanofin}[2]{\spanof{#1}^{#2}}
We now deduce Theorems \ref{thm:categorical0} and \ref{thm:categorical} by
proving quasi-minimal excellence of an appropriate class of structures
corresponding to the exact sequences of (\ref{exactSeqChar0}) and
(\ref{exactSeqCharp}). We use \cite{KirbyQME} as our reference for the theory
of quasi-minimal excellent classes.

The following is a corrected and abbreviated version of the argument in
\cite[Section 3]{ZCovers}.

Let $L$ be the one-sorted language
$\seq{+,(\mu_q)_{q\in\Rats},\omega,(W_f)_{(f\in\Ints[X_1,\ldots,X_n],
n\in\Nats)},E}$.
If $p>0$, fix a map $\ex_0$ as in (\ref{eqn:ex_0}).  Let $\Sigma$ be the
$L_{\omega_1,\omega}(L)$-sentence expressing that for a model $V$:
\begin{enumerate}[(I)]
  \item $(V;+,(\mu_q)_{q\in\Rats})$ is a $\Rats$-vector space (we write $qx$
    for $\mu_q(x)$);
  \item $E$ is an equivalence relation on $V$;
  \item $V/E$ can be identified with the multiplicative group
    $\multgrpof{F}$ of a characteristic $p$ algebraically closed field
    $(F;+,\cdot)$ such that $(x+y)/E = x/E \cdot y/E$, and for each
    $n\in\Nats$ and each polynomial
    $f \in \Ints[X_1,\ldots,X_n]$, we have $V\models W_f(x_1,\ldots,x_n)$ iff
    $f(x_1/E,\ldots,x_n/E) = 0$.
  \item $xE0$ iff $\bigvee_{z\in\emp} x = z\omega$.
  \item (if $p>0$) for each tuple of rationals $\qtup\in\Rats$, the field
    types of $(q_1\omega/E,\ldots,q_n\omega/E)$ and
    $(\ex_0(q_1),\ldots,\ex_0(q_n))$ are equal.
\end{enumerate}

Models of $\Sigma$ correspond to exact sequences as in (\ref{exactSeqChar0})
and (\ref{exactSeqCharp}); given a model $V\models\Sigma$ we write $\ex :
V \maps \multgrpof{F}$ for the quotient map. For $p>0$, axiom (V) implies that
for an appropriate choice of embedding $\Fpalg\leq\multgrpof{F}$, we have that
$\ex$ extends $\ex_0$.

Let $\CClass$ be the class of models of $\Sigma$.

We equip $V\in\CClass$ with a closure operation $\cl(X) :=
\ex^{-1}(\acl(\ex(X)))$; as in \cite[Lemma 3.2]{ZCovers}, 
this satisfies \cite[Axiom I]{KirbyQME} as well as the exchange and
countable closure properties.

In the following, a \defnstyle{partial embedding} is a partial map $f$ which
extends to an isomorphism $\spanof{f} : \spanof{\dom(f)} \maps
\spanof{\image(f)}$, where for $A\subseteq V\in\CClass$, we define $\spanof{A} =
\spanofin{A}{V}$ to be the substructure of $V$ generated by $A$. By our
choice of language, $\spanof{A}$ is the $\Rats$-vector space span of
$A\union\{\omega\}$.

Although it is not specified by the axioms in \cite{KirbyQME}, the following
property is in fact necessary for the categoricity theorems \cite[Theorem 2.1,
Theorem 3.3]{KirbyQME}:
\begin{lemma}
  If $V_1,V_2\in\CClass$ then the substructures generated by the empty set,
  $\spanofin{\emptyset}{V_i}\leq V_i$, are isomorphic.
\end{lemma}
\begin{proof}
  In positive characteristic, this is immediate from Axiom (V).

  In characteristic $0$, we argue as follows. The map
  \begin{align*}
    \theta: \mu_1 &\maps \mu_2\\
    \theta(\ex_1(q\omega_1)) &:= \ex_2(q\omega_2)
  \end{align*}
  is a group isomorphism of the torsion groups. It follows (see \cite[VI
  3.1]{LangAlg}) that $\theta$ is a partial field isomorphism; hence
  $q\omega_1 \mapsto q\omega_2$ is an isomorphism
  $\spanofin{\emptyset}{V_1} \maps \spanofin{\emptyset}{V_2}$ as required.
\end{proof}

We proceed to verify Axioms II and III of \cite{KirbyQME}.

The following lemma proves \cite[Axiom II]{KirbyQME}.

\begin{lemma}[$\omega$-homogeneity over submodels and $\emptyset$]
  Let $V_1,V_2\in \CClass$, let $G_i\subseteq V_i$ be closed
  substructures or the empty set, and let $g : G_1 \maps G_2$ be an isomorphism or
  the empty map.
  \begin{enumerate}[(i)]
    \item If $x_i \in V_i \setminus \cl(G_i)$, then $g\union\{(x_1,x_2)\}$ is a
      partial embedding.
    \item If $\atup_1\in V_1$ and $g': G_1\atup_1\maps V_2$ is a partial
      embedding extending $g$, then for any $b_1\in\cl(G_1\atup_1)$ there
      exists $b_2\in V_2$ such that $g'\union\{(b_1,b_2)\}$ is a partial
      embedding.
  \end{enumerate}
\end{lemma}
\begin{proof}
  We have $\ex_i : V_i \maps \multgrpof{F_i}$.

  (i) is clear.

  For (ii), suppose first that $G_i$ is closed, so $\ex_i(G_i) =
  \multgrpof{F_i'}$ where $F_i'\leq F_i$ is an algebraically closed subfield.
  
  We may assume that $(\atup_1b_1)$ is linearly independent over $G_1$.
  By the $n=1$ case of Theorem \ref{thm:determined}, division systems below
  $\ex(b_1)$ are finitely determined over $\ex_1(\spanof{G_1,\atup_1})$. The
  result follows.
  
  The case remains that $G_i=\emptyset$. In characteristic 0, we refer
  to \cite[3.5(ii)]{ZCovers} for this. In characteristic $p>0$, the
  substructure of $G_i$ generated by $\emptyset$ is a closed subset of $V_i$,
  and so we return to the case above.
\end{proof}

The following lemma proves \cite[Axiom III]{KirbyQME} for $\CClass$ - that
axiom refers to {\em closed} partial embeddings, but in $\CClass$ any partial
embedding is closed.

\begin{lemma}
  Let $V_1,V_2\in\CClass$, let $B\subseteq V_1$ be a $\cl$-independent set, let
  $B_1,\ldots,B_n\subseteq B$, and let $C := \bigcup \cl(B_i) \subseteq V_1$.
  Let $g:C \maps V_2$ be a partial embedding.
  Let $\atup\in\cl(C)$. Then there exists a finite subset $C_0$ of $C$ such
  that if $g' : C_0\atup \maps V_2$ is a partial embedding extending $g\restriction_{C_0}$
  then $g\union g' : C\atup \maps V_2$ is a partial embedding.
\end{lemma}
\begin{proof}
  Letting $L_i := \ex(\cl(B_i))$, we have that $L_1,\ldots,L_n$ are from an
  independent system. We may assume that $\atup$ is $\Rats$-linearly
  independent over $\Sigma_i \cl(B_i)$. By Theorem \ref{thm:determined}, division
  systems below $\ex(\atup)$ are finitely determined over $\bigvee_i L_i$.
  Let $m$ be as in the definition of finite determination; $\ex(\atup/m)$ is
  algebraic over $\bigvee_i L_i$, and so its field type is isolated by some
  field formula $\phi(\xtup,\btup_1,\ldots,\btup_n)$ where $\btup_i\in L_i$.
  Letting $C_0\subseteq C$ be a finite subset such that $\btup_i\in\ex(C_0)$
  for all $i$, we see that $C_0$ is as required.
\end{proof}

We have now shown that $\CClass$ is a quasiminimal excellent class in the
sense of $\cite{KirbyQME}$. By \cite[Theorem 3.3]{KirbyQME}, therefore, there is at
most one structure in $\CClass$ of a given $\cl$-dimension, i.e. with the
corresponding algebraically closed field having a given transcendence degree.
That there exists such a structure in each transcendence degree is clear.
Translating straightforwardly from our one-sorted setup to the two-sorted
setup of their statements, this concludes the proofs of Theorems
\ref{thm:categorical0} and \ref{thm:categorical}.

\bibliography{mt,myBooks,misc}

\end{document}